%% file: proof.tex
\documentclass{amsart}
\usepackage{amssymb}

\newtheorem{thm}{Theorem}
\newtheorem{prop}{Proposition}

\def\bbR{\mathbb{R}}

\newcommand\pd[2]{\displaystyle\frac{\partial #1}{\partial #2}}

\DeclareMathOperator{\lk}{lk}

\begin{document}
\title{Linking Integral Projection}
\author{Daniel J. Cross}
\address{Drexel University}
\email{d.j.cross@drexel.edu}
\begin{abstract}
The linking integral is an invariant of the link-type of two manifolds immersed in a Euclidean space.  It is shown that the ordinary Gauss integral in three dimensions may be simplified to a winding number integral in two dimensions.  This result is then generalized to show that in certain circumstances the linking integral between arbitrary manifolds may be similarly reduced to a lower dimensional integral. 
\end{abstract}

\maketitle

\section{Reduction of the Gauss Integral to the Winding Number Integral}
\label{sec:wind}
The linking number of two disjoint oriented closed curves in $\bbR^3$ is an
integer invariant that in some sense measures the extent of linking
between the curves.  While there are many equivalent ways to compute
this number\cite{Rol03}, the most well-known is the linking integral of Gauss.  In
this section we show that this integral in 3-space may always be
simplified to an integral in 2-space which is equivalent to a winding
number integral.  

\begin{prop}Given two disjoint immersed closed curves $s\mapsto \gamma_1(s)$ and
$t\mapsto\gamma_2(t)$ in $\bbR^3$, the Gauss linking integral of the pair reduces
to a sum of winding numbers of one curve about a sequence of points
determined by the other, contained in some 2-dimensional hyperplane.
\end{prop}
\begin{proof}
The link of $\gamma_1$ and $\gamma_2$, $\lk(\gamma_1,\gamma_2)$, is given by the Gauss integral,
$$
\lk(\gamma_1,\gamma_2)=\frac{1}{4\pi}\int\det\left(r,\pd{\gamma_1}{s},\pd{\gamma_2}{t}\right)\frac{dsdt}{r^3},
$$
where each term in the determinant is a column vector and
$r=\gamma_2-\gamma_1$ is the relative position vector.
Through a homotopy of the maps we may arrange $\gamma_1$ to lie in the
plane $x_3=0$ with $\gamma_2$ intersecting the plane perpendicularly in a
finite number of points $p_i$.  This may be done so that $\gamma_1$ and
$\gamma_2$ remain disjoint throughout.  Since the Gauss integral is a
homotopy invariant, $\lk(\gamma_1,\gamma_2)$ is preserved through this
deformation.  If the homotopy was merely continuous we may replace
it with an arbitrarily close smooth homotopic approximation.

Next, deform $\gamma_2$ near the intersection with the plane so that
it becomes a straight line segment perpendicular to the plane in a
neighborhood of each intersection point.  Now deform it further by
``stretching'' it away from the plane so that the straight line
segments are extended further away from the plane and the rest of
$\gamma_2$ is pushed further away from the plane.  In the limit that
the stretching goes off to infinity, the denominator of the integral
falls off sufficiently fast that its contribution goes to zero.  We
are left with a finite number of infinite line segments perpendicular
to the plane and disjoint from $\gamma_1$. We assume each line is
parameterized in the standard way, $t\mapsto \pm x_3$.

We will now assume each line parameterized by $t\mapsto x_3$ so that
$\partial\gamma_2/\partial t=e_3$, but introduce an orientation to each
point $o(p_i)$ which is $\pm 1$ depending on the original
parameterization of the corresponding line in an obvious way.  We then
see that the linking integral becomes
$$
\sum_io(p_i)\frac{1}{4\pi}\int\det\left(r,\pd{\gamma_2}{s},e_3\right)\frac{dsdt}{(\rho^2+t^2)^{3/2}},
$$
where $\rho$ is the restriction of $r$ to the plane $x_3=0$.

Notice that since $e_3=(0,0,1)^t$ the
determinant reduces to that of the upper-left block, which is
$\det(\rho_{12},\partial\gamma_2/\partial s))$ and is independent of
$t$.  Thus we may evaluate the integral
$$
\int_\bbR\frac{dt}{(\rho^2+t^2)^{3/2}}=\frac{2}{\rho^2},
$$
and the linking integral becomes
$$
\sum_io(p_i)\frac{1}{2\pi}\int\det\left(\rho,\pd{\gamma_2}{s}\right)\frac{ds}{\rho^2},
$$
which is easily seen to be the sum of the winding numbers of
$\gamma_2$ about each point $p_i$ times the orientation of $p_i$.
\end{proof}

The construction in the proof also allows one to show the linking
integral may also be given as an intersection number of $\gamma_1$
with a surface spanned by $\gamma_2$.  Indeed, perturb $\gamma_2$ to
an embedding and let $S$ be a Seifert surface constructed by Seifert's
algorithm\cite{Seifert34,Rol03}.  The number of Seifert discs above a $p_i$ is precisely 
the winding number of $\gamma_1$ about $p_i$ and the induced
orientation of each Seifert disc is given by the orientation of the
bounding curve.  Finally, with $o(p_i)$ we have the signed
intersection number of $\gamma_2$ with the Seifert disc, and the sum
over all gives the signed intersection number of $\gamma_2$ with $S$.

\section{The General Linking Integral Projection}
In this section the proposition of section \ref{sec:wind}
is generalized
from curves to arbitrary compact boundaryless oriented manifolds $M^n$ and
$N^n$ mapped disjointly into $\bbR^{p+1}$, $p=m+n$.  In this case one may define a
linking number by  ${\lk}(M,N)=(-1)^{m}\deg \hat r$, where $\hat r$
is the unit relative position vector defined by 
\begin{align*}
\hat r:M\times N&\to S^p\\
(x,y)&\mapsto \frac{r}{||r||}=\frac{x-y}{||x-y||},
\end{align*}
where $x$ and $y$ are points in the images of $M$ and $N$ in
$\bbR^{p+1}$ respectively.  We will show that under certain conditions
the linking number calculation reduces to a calculation in a
hyperplane.  We note that these expressions may be defined with
differing sign conventions in which case the conclusion of the theorem
will hold up to a sign. The present convention is most convenient for
expressing the present result.

\begin{thm}
Given $M$ and $N$ as above, suppose that there exists smooth homotopies of
$M$ and $N$  maintaining disjointness and taking $M$ into an
$m+n'+1$-dimensional  hyperplane $H$, $0\leq n'\leq n$, and that $H$ 
intersects $N$ transversely in the submanifold $N'$.  Then
$\lk(M,N)=\lk(M,N')$, where the first linking integral is taken in
$\bbR^{p+1}$ and the second in $H\cong\bbR^{p'+1}$, where $p=m+n$ and $p'=m+n'$. 
\end{thm}
\begin{proof}
It is straightforward to show\,\cite{frankel06} that the degree of this map may
be written explicitly as
\begin{equation}\label{eq:app_link}
\deg\hat r = \frac{(-1)^m}{{\rm vol}S^p}\int_{M\times N}
\det\left(r,\pd{x}{s},\pd{y}{t}\right)
\frac{dsdt}{||r||^{p+1}},
\end{equation}
where $s$ and $t$ represent oriented local coordinates $s_i$ and $t_j$
on $M$ and $N$ respectively, and the quantities in the determinant are
column vectors. 

We now homotope $M$ into $H\simeq\bbR^{p'+1}$ and
homotope $N$ so that it intersects $H$ transversely.
The intersection $N'=\cup N'_i$ will be a finite disjoint union of
closed oriented manifolds  of dimension $n'$ (the codimension of the
transverse intersection of two manifolds is the sum of their
codimensions).  We may
actually assume that $N\perp H$ (in the Euclidean metric of
$\bbR^{p+1}$) so that $N$ is locally of the form $N'\times\bbR^{n-n'}$ in
some neighborhood of $H$.  Now extend this local product
decomposition by pushing the rest of $N$ off to $\infty$ as was done
in section \ref{sec:wind}.  Since the 
integrand in Eq.\,\ref{eq:app_link} falls off sufficiently fast with distance, this contribution
to the integral goes to zero, so we may make the replacement $N\to
N'\times\bbR^{n-n'}$.

Adapt the coordinates on $N$ with respect to the product
decomposition so that the last $n-n'$ coordinates are Euclidean
coordinates on $\bbR^{n-n'}$.  The partial derivatives $\partial
y/\partial t_i$ with respect to these coordinates are just $\pm e_i$, the
$i$th unit vector, but the signs may vary on different $N'_i$.  We
may absorb these signs into the orientation of the components,
considering their orientations reversed if necessary (rather than
explicitly introducing an orientation function as was done in section
\ref{sec:wind}).

The matrix in the integrand has a block
structure with an $n-n'$ unit matrix in the lower right block and zero
in the upper right block.  Thus the determinant may be replaced with
just that of the upper left block.  Since this matrix is independent of the
last $n-n'$ coordinates the distance function $r$ reduces to $\rho=r|_H$.
It remains to evaluate the integral
\begin{equation*}
I=\int_{-\infty}^\infty \cdots\int_{-\infty}^\infty \frac{dt_{q+1}\cdots dt_n}{\left(\rho^2+\sum t_i^2\right)^{(p+1)/2}},
\end{equation*}
for $i=q+1,\ldots,n$.  Write $t_n=z$ and $a=\rho^2+\sum t_i^2$, $i\neq
n$ and then
\begin{equation*}
\int_{-\infty}^\infty \frac{dz}{\left(a+z^2\right)^{(p+1)/2}}=
\frac{\sqrt{\pi}}{a^{p/2}}\frac{\Gamma(\frac{p}{2})}{\Gamma(\frac{p+1}{2})}
=\frac{1}{a^{p/2}}\frac{{\rm vol}_{S^p}}{{\rm   vol}_{S^{p-1}}},
\end{equation*}
using the well-known expression for the volume of a sphere.  By
progressively isolating the variables $t_i$ we obtain an 
integral of the same form but with $p$ decreasing by one each time.
Proceeding by induction we obtain 
\begin{equation*}
I=\frac{1}{||\rho||^{p'+1}}\frac{{\rm vol}_{S^p}}{{\rm vol}_{S^{p'}}},
\end{equation*}
where $p'=p-n+n'=m+n'$.  Hence Eq.\,\ref{eq:app_link} becomes 
\begin{align*}
\deg\hat r &= \frac{(-1)^m}{{\rm vol}S^{p'}}\int_{M\times N'}
\det\left(\rho,\pd{x}{t},\pd{y}{s}\right)
\frac{dsdt}{||\rho||^{p'+1}}\\
&=(-1)^{m}\deg\hat\rho,
\end{align*}
which is $\lk(M,N')$.
\end{proof}

This result has been applied in \cite{Cross09a} to demonstrate that two
solid tori ($D^2\times S^1$) embedded into $\bbR^4$ differing by a
Dehn twist are non-isotopic.  In fact, it is shown that the isotopy
classes of embeddings of the solid torus into $\bbR^4$ are in
bijective correspondence with the integers, the correspondence given
by the number of applied Dehn twists.

\input{proof.bbl}

\end{document}

%% file: proof.bbl
\providecommand{\bysame}{\leavevmode\hbox to3em{\hrulefill}\thinspace}
\providecommand{\MR}{\relax\ifhmode\unskip\space\fi MR }
\providecommand{\MRhref}[2]{%
  \href{http://www.ams.org/mathscinet-getitem?mr=#1}{#2}
}
\providecommand{\href}[2]{#2}